\documentclass[12pt]{amsart}
\usepackage{amssymb}
\usepackage[OT2,T1]{fontenc}
\DeclareSymbolFont{cyrletters}{OT2}{wncyr}{m}{n}
\DeclareMathSymbol{\Sha}{\mathalpha}{cyrletters}{"58}
\usepackage{amssymb,amsmath,amsthm, hyperref, verbatim, amsfonts, fontenc}
\newtheorem*{thm*}{Theorem}
\newtheorem{thm}{Theorem}[section]
\newtheorem{theorem}{Theorem}[section]

\newtheorem*{lemma}{Lemma}
\numberwithin{equation}{section}
\newtheorem{prop}{Proposition}[section]
\newtheorem*{prop*}{Proposition}

\theoremstyle{remark}
\newtheorem*{rmk}{Remark}

\theoremstyle{definition}

\newtheorem*{dfn}{Definition}

\newtheorem{corollary}[thm]{Corollary}
\newtheorem*{corollary*}{Corollary}

\newcommand{\C}{\mathbb{C}}
\newcommand{\Q}{\mathbb{Q}}

\newcommand{\Z}{\mathbb{Z}}

\newcommand{\R}{\mathbb{R}}

\newcommand{\Comment}[1]{}

\begin{document}
\title{A ``Strange'' Vector-Valued Quantum Modular Form}
\author{Larry Rolen and Robert P. Schneider}
\email{larryrolen@gmail.com}
\email{robert.schneider@emory.edu}
\address{Department of Mathematics and Computer Science\newline
Emory University\newline
400 Dowman Dr., W401\newline
Atlanta, GA 30322}
\thanks{The authors would like to thank our advisor Ken Ono for his guidance throughout this project. We are also grateful to Robert Rhoades and to the referee for their helpful comments; in particular the referee suggested using Fine's identity for an alternative expression for $\theta_1^S(q^{-1})$ in (\ref{inversefine}).}
\maketitle
\begin{abstract}
Since their definition in 2010 by Zagier, quantum modular forms have been connected to numerous different topics such as strongly unimodal sequences, ranks, cranks, and  asymptotics for mock theta functions near roots of unity. These are functions that are not necessarily defined on the upper half plane but \emph{a priori} are defined only on a subset of $\Q$, and whose obstruction to modularity is some analytically ``nice'' function. Motivated by Zagier's example of the quantum modularity of Kontsevich's ``strange'' function $F(q)$, we revisit work of Andrews, Jim{\'e}nez-Urroz, and Ono to construct a natural vector-valued quantum modular form whose components are similarly ``strange''. 
 \end{abstract}

\section{Introduction and Statement of Results}
In his seminal 2010 Clay lecture, Zagier defined a new class of function with certain automorphic properties called a ``quantum modular form'' \cite{Zagier-Clay}. Roughly speaking, this is a complex function on the rational numbers which has modular transformations modulo ``nice'' functions. Although the definition is intentionally vague, Zagier gave a handful of motivating examples to serve as prototypes of quantum behavior. For example, he defined quantum modular forms related to Dedekind sums, sums over quadratic polynomials, Eichler integrals, and to Kontsevich's ``strange'' function. More concisely, Zagier proposed the following:
\begin{dfn} We say that a function $f$ from a subset of $\mathbb{P}^1(\Q)$ to $\C$ is a \emph{quantum modular form} if $f(x)-f|_k\gamma(x)=h_{\gamma}(x)$ for a ``suitably nice'' function $h_{\gamma}(x)$.
\end{dfn}
Here $|_k$ is the usual Petersson slash operator, and ``suitably nice'' implies some pertinent analyticity condition, e.g. $\mathcal{C}_k,\mathcal{C}_{\infty}$, etc. One of the most striking examples of quantum modularity is described in Zagier's paper on Vassiliev invariants \cite{Zagier-Vassiliev}, in which he studies the Kontsevich ``strange'' function $F(q)$ given by

\begin{equation}F(q):=\sum_{n=0}^{\infty}(q;q)_n,\end{equation}
where $(a,t)_n:=\prod_{i=0}^{n-1}(1-at^i)$ denotes the usual \emph{$q$-Pochhammer symbol} with the conventions $(a,t)_0:=1$, $(a;t)_{\infty}:=\lim_{n\rightarrow\infty}(a;t)_n$, and $q:=e^{2\pi i z}$ with $z\in\C$.\\ \indent This function is strange indeed, as it does not converge on any open subset of $\C$, but converges as a finite sum for $q$ any root of unity. Recently, Bryson, Pitman, Ono, and Rhoades showed \cite{Unimodal} that $F(q^{-1})$ agrees to infinite order at roots of unity with a function $U(-1,q)$ which is also well-defined on the upper-half plane $\mathbb{H}$, obtaining a quantum modular form that is a ``reflection'' of $F(q)$ and that naturally extends to $\mathbb{H}$. Moreover, $U(-1,q)$ counts unimodal sequences having a certain combinatorial statistic. 
\\

Zagier's study of $F(q)$ depends on the formal $q$-series identity
\begin{equation}
\displaystyle\sum_{n=0}^{\infty}\left(\eta(24z)-q(1-q^{24})(1-q^{48})\cdots(1-q^{24n})\right)=\eta(24z)D(q)+E(q),
\end{equation}
\noindent
where $\eta(z):=q^{1/24}(q;q)_{\infty}$, $D(q)$ is an Eisenstein-type series, and $E(q)$ is a ``half-derivative'' of $\eta(24z)$ (such formal half-derivatives will be discussed in $\S2$). We refer to such an identity as a ``sum of tails'' identity. 
Here we revisit Zagier's construction using work of Andrews, Jim{\'e}nez-Urroz, and Ono on more general sums of tails formulae \cite{A-J-O} (see also \cite{Andrews-Advances}). We construct a natural $3$-dimensional vector-valued quantum modular form associated to tails of infinite products. Moreover, the components are analogous ``strange'' functions; they do not converge on any open subset of $\C$ but make sense for an infinite subset of $\Q$. We define:
\begin{equation}H(q)=\begin{pmatrix}\theta_1\\ \theta_2 \\ \theta_3\end{pmatrix}:=\begin{pmatrix}\eta(z)^2/\eta(2z)\\ \eta(z)^2/\eta(z/2)\\ \eta(z)^2/\eta(\frac{z}{2}+\frac12) \end{pmatrix}. \end{equation}
We also note that $\theta_3=\zeta_{48}^{-1}\cdot\frac{\eta(z/2)\eta(2z)}{\eta(z)}$ by the following identity which is easily seen by expanding the product definition of $\eta(z)$:
\begin{equation}\label{eta_identity1}
\eta(z+1/2)=\zeta_{48}\cdot\frac{\eta(2z)^3}{\eta(z)\cdot\eta(4z)},
\end{equation}
where $\zeta_k:=e^{2\pi i/k}$. From this it follows that if we let 
\begin{equation*}
\quad\quad F_9(z):=\eta(z)^2/\eta(2z),\quad\quad\quad\quad\quad\quad\quad\quad
F_{10}(z):=\eta(16z)^2/\eta(8z)
\end{equation*}
then 
\begin{equation}\label{F9F10}
H(q)=\begin{pmatrix}F_9(q)&F_{10}(q^{1/16})
&\zeta_{12}^{-1}F_{10}(\zeta_{16}\cdot q^{1/16})\end{pmatrix}^T\end{equation}
(the notations $F_9$ and $F_{10}$ come from \cite{A-J-O}).
For convenience, we recall the classical theta-series identities for $F_9$ and $F_{10}$:
\begin{equation}
F_9(q)= 1+2\displaystyle\sum_{n=1}^{\infty}(-1)^nq^{n^2},     \quad\quad\quad\quad\quad\quad\quad\quad F_{10}(q)=\displaystyle\sum_{n=0}^{\infty}q^{(2n+1)^2}
.\end{equation}
It is simple to check that $H(z)$ is a $3$-dimensional vector-valued modular form using basic properties of $\eta(z)$, as we describe in $\S4$. To each component $\theta_i$ we associate for all $n\geq0$ a finite product $\theta_{i,n}$:
\begin{equation}
\theta_{1,n}:=\frac{(q;q)_n}{(-q;q)_n},\quad\quad\quad \theta_{2,n}:=q^{\frac{1}{16}}\cdot\frac{(q;q)_n}{(q^{\frac12};q)_{n+1}},\quad\quad\quad\theta_{3,n}:=\frac{\zeta_{16}}{\zeta_{12}}\cdot q^{\frac{1}{16}}\cdot\frac{(q;q)_n}{(-q^{\frac12};q)_{n+1}},
\end{equation}
such that $\theta_{i,n}\rightarrow\theta_i$ as $n\rightarrow\infty$.
Next, we construct corresponding ``strange'' functions
$
\theta_i^S:=\sum_{n=0}^{\infty}\theta_{i,n}.
$
Note that these functions do not make sense on any open subset of $\C$, but that each $\theta_i^S$ is defined for an infinite set of roots of unity and, in particular, $\theta_2^S$ is defined for all roots of unity.
Our primary object of study will then be the vector of ``strange'' series $H_Q(z):=\begin{pmatrix}\theta_1^S(z)& \theta_2^S(z)& \theta_3^S(z)\end{pmatrix}^T.$ In order to obtain a quantum modular form, we first define $\phi_i(x):=\theta_i^S(e^{2 \pi i x})$ from a subset of $\Q$ to $\C$, and let $\phi(x):=\begin{pmatrix}\phi_1(x)& \phi_2(x)& \phi_3(x)\end{pmatrix}^T.$ We then show the following result.
\begin{theorem}\label{mainthm} Assume the notation above. Then the following are true:
\begin{enumerate} 
\item There exist $q$-series $G_i$ (see $\S4$) which are well-defined for $|q|<1$, such that $\theta_i^S(q^{-1})=G_i(q)$ for any root of unity for which $\theta_i^S$ converges.
\item We have that $\phi(x)$ is a weight $3/2$ vector-valued quantum modular form. In particular, we have that 
\[
\phi(z+1)-\begin{pmatrix}1&0&0\\0&0&\zeta_{12}\\0&\zeta_{24}&0 \end{pmatrix}\phi(z)=0
,
\]
and we also have that
\[
\left(\frac{z}{-i}\right)^{-3/2}\phi(-1/z)+\begin{pmatrix}0&\sqrt{2}&0\\1/\sqrt{2}&0&0\\0&0&1\end{pmatrix}\phi(z)=\begin{pmatrix}0&\sqrt{2}&0\\1/\sqrt{2}&0&0\\0&0&1\end{pmatrix}g(z)
,
\]
where $g(z)$ is a $3$-dimensional vector of smooth functions defined as period integrals in $\S3$.
\end{enumerate}
\end{theorem}
In addition, we deduce the following corollary regarding generating functions of special values of zeta functions from the sums of tails identities. Let
\begin{equation}
H_9(t,\zeta):=-\frac14\displaystyle\sum_{n=0}^{\infty}\frac{(1-\zeta e^{-t})(1-\zeta^2e^{-2t})\cdots(1-\zeta^n e^{-nt})}{(1+\zeta e^{-t})(1+\zeta^2e^{-2t})\cdots(1+\zeta^ne^{-nt})},
\end{equation}
\begin{equation}
H_{10}(t,\zeta):= -2(\zeta e^{-t})^{1/8}\displaystyle\sum_{n=0}^{\infty}\frac{(1-\zeta e^{-2t})(1-\zeta^2e^{-4t})\cdots(1-\zeta^n e^{-2nt})}{(1-\zeta e^{-t})(1-\zeta^2e^{-3t})\cdots(1-\zeta^ne^{-(2n+1)t})}.
\end{equation}
\begin{rmk}
Note that there are no rational numbers for which all three components of $\phi$ make sense simultaneously. To be specific, $\phi_1(z)$ makes sense for rational numbers which correspond to primitive odd order roots of unity, $\phi_2(z)$ makes sense for all rational numbers, and $\phi_3(z)$ converges at even order roots of unity. Hence, by the equation in (2) of \ref{mainthm}, we understand that each of the six equations of the vector-valued transformation laws is true where the corresponding component in the equation is well-defined; as there are no equations in which $\phi_1$ and $\phi_3$ both appear, then for all the equations there is an infinite subset of rationals on which this is possible.
\end{rmk}

For a root of unity $\zeta$, we define the following two $L$-functions

\[L_1(s,\zeta):=\sum_{n=1}^{\infty}\frac{(-\zeta)^{n^2}}{n^s},\]

\[L_2(s,\zeta):=\sum_{n=1}^{\infty}\left(\frac2n\right)^2\cdot\frac{\zeta^{\frac{n^2}{8}}}{n^s}.\]
Then we have the following.

\begin{corollary}
Let $\zeta=e^{2\pi i \alpha}$ be a primitive $k^{\operatorname{th}}$ root of unity, $k$ odd for $H_9$ and $k$ even for $H_{10}$. Then as $t\searrow0$, we have as power series in $t$

\begin{equation}\label{Hurwitz1}  H_{9}(t,\zeta)=\sum_{n=0}^{\infty}\frac{L_1(-2n-1,\zeta)(-t)^n}{n!},
\end{equation}
\begin{equation}\label{Hurwitz2} H_{10}(t,\zeta)=\sum_{n=0}^{\infty}\frac{L_2(-2n-1,\zeta)(-t)^n}{8^nn!}.
\end{equation}
\end{corollary}

To illustrate our results by way of an application, we provide a  numerical example which gives finite evaluations of seemingly complicated period integrals. First define 
\begin{equation*}
\Omega(x):=\int_{x}^{i\infty}\frac{\theta_1(z)}{(z-x)^{3/2}}\, dz
\end{equation*} 
for $x\in\Q$, and consider $\theta_1^S(\zeta_k)$ for $k$ odd, which is a finite sum of $k^{\text{th}}$ roots of unity. Then the proof of Theorem \ref{mainthm} will imply that $\Omega(1/k)=\pi i(1+i)\theta_1^S(\zeta_k)$ by showing that the period integral $\Omega(x)$ is a ``half-derivative'' which is related to $\theta_1^S$ at roots of unity by a sum of tails formula. The following table gives finite evaluations of $\theta_1^S(\zeta_k)$ and numerical approximations to the integrals $\Omega(1/k)$:\\ 
\begin{center} 
\begin{tabular}{c c c c} 
\hline\hline                        
$k$&$\pi i(i+1)\theta_1^S(\zeta_k)$&$\int_{1/k+10^{-9}}^{10^9i}\frac{\theta_1(z)}{(z-1/k)^{\frac32}}\, dz$\\ [0.7ex] 
\hline                    
\\
$3$&$\pi i(i+1)(-2 \zeta_{3} + 3)\sim-7.1250+18.0078i$&$-7.1249+18.0078i$
\\
$5$&$\pi i(i+1)(-2 \zeta_{5}^{3} - 2 \zeta_{5}^{2} - 8 \zeta_{5} + 3)\sim12.078+35.7274i$&$12.078+35.7273i$
\\
$7$&$\pi i(i+1)(6 \zeta_{7}^{4} - 2 \zeta_{7}^{2} - 10 \zeta_{7} + 7)\sim52.0472+25.685i$&$52.0474+25.685i$
\\
$9$&$\pi i(i+1)(8 \zeta_{9}^{4} - 16 \zeta_{9} + 3)\sim76.4120-28.9837i$&$76.4116-28.9836i$
\\
\hline     
\end{tabular}
\end{center}
\vspace{0.2in}
The paper is organized as follows. In $\S2$ we recall the identities of \cite{A-J-O}, and in $\S3$ we describe the modularity properties of Eichler integrals of half-integral weight modular forms. In $\S4$ we complete the proof of Theorem \ref{mainthm}. We finish with the proof of Corollary $1.1$ in $\S5$.
\section{Preliminaries}
In this section, we describe some of the machinery needed to prove Theorem \ref{mainthm}
\subsection{Sums of Tails Identities}
Here we recall the work of Andrews, Jim{\'e}nez-Urroz, and Ono on sums of tails identities. To state their results for $F_9$ and $F_{10}$ and connect $\theta_i^S$ to quantum modular objects, we formally define a ``half-derivative operator'' by
\begin{equation}
\sqrt{\theta}\left(\displaystyle\sum_{n=0}^{\infty}a(n)q^n\right):=\displaystyle\sum_{n=1}^{\infty}\sqrt{n}a(n)q^n.
\end{equation}
If we have a generic $q$-series $f(q)$, we will also denote $\sqrt{\theta}f(q):=\widetilde{f}(q)$.
Then Andrews, Jim{\'e}nez-Urroz, and Ono show \cite{A-J-O} that for finite versions $F_{9,i}, F_{10,i}$ associated to $F_{9},F_{10}$ the following holds true:
\begin{theorem}[Andrews-Jim{\'e}nez-Urroz-Ono]As formal power series, we have that
\begin{equation}
\displaystyle\sum_{n=0}^{\infty}\left(F_{9}(z)-F_{9,n}(z)\right)=2F_{9}(z)E_1(z)+2\sqrt{\theta}(F_{9}(z))
,\end{equation}

\begin{equation}
\displaystyle\sum_{n=0}^{\infty}\left(F_{10}(z)-F_{10,n}(z)\right)=F_{10}(z)E_2(z)+\frac12\sqrt{\theta}(F_{10}(z)),
\end{equation}
\end{theorem} 
\noindent
where the $E_i(z)$ are holomorphic Eisenstein-type series. In particular, as $F_9,F_{10}$ vanish to infinite order while $E_1,E_2$ are holomorphic at all cusps where the ``strange'' functions are well-defined, we have for $q$ an appropriate root of unity that the ``strange'' function associated to $F_i$ equals $\widetilde{F_i}$ to infinite order. As the series $\theta_2,\theta_3$ do not have integral coefficients, we make the definitions $\widetilde{\theta}_2(z):=\widetilde{F}_{10}(z/16)$ and $\widetilde{\theta}_3(z):=\widetilde{F}_{10}(z/16+1/16)$. By the definition of the strange series, we obtain the following. 
\begin{corollary}
At appropriate roots of unity where each ``strange'' series is defined, we have that
\begin{equation}
\theta_1^S(q)=2\widetilde{\theta_1}(q),\quad\quad\quad\theta_2^S(q)=\frac12\widetilde{\theta_2}(q),\quad\quad\quad\theta_3^S(q)=\frac12\widetilde{\theta_3}(q).
\end{equation}
\end{corollary}
\section{Properties of Eichler Integrals}
In the previous section we have seen that at a rational point $x$, each component of $\phi(x)$ agrees up to a constant with a ``half-derivative'' of the corresponding theta function at $q=e^{2\pi i x}$. Thus, we can reduce part (2) of Theorem \ref{mainthm} to a study of modularity of such half-derivatives. We do so following the outline given in \cite{Zagier-Vassiliev}, which is further explained in the weight $3/2$ case in \cite{L-Zagier}. Recall that in the classical setting of weight $2k$ cusp forms, $1\leq k\in\Z$, we define the \emph{Eichler integral} of $f(z)=\sum_{n=1}^{\infty} a(n)q^n$ as a formal $(k-1)^{\text{st}}$ antiderivative
$
\widetilde{f}(z):=\sum_{n=1}^{\infty} n^{1-k}a(n)q^n.
$
Then $\widetilde{f}$ is nearly modular of weight $2-k$, as the differentiation operator $\frac{\, d}{\, dq}$ does not preserve modularity but preserves near-modularity. More specifically, $\widetilde{f}(z+1)=\widetilde{f}(z)$ and $z^{k-2}\widetilde{f}(-1/z)-\widetilde{f}(z)=g(z)$ where $g(z)$ is the \emph{period polynomial}. This polynomial encodes deep analytic information about $f$ and can also be written as $g(x)=c_k\int_0^{i\infty}f(z)(z-x)^{k-2}\, dz$ for a constant $c_k$ depending on $k$. Suppose we now begin with a weight $1/2$ vector-valued modular form $f$ with $n$ components $f_i$ such that and $f(-1/z)=M_Sf(z)$, for $M_S$ both $n\times n$ matrices (the transformation under translation is routine).  
\\
\indent
In this case, of course, it does not make sense to speak of a half-integral degree polynomial, and the integral above does not even converge. However, we may remedy the situation so that the analysis becomes similar to the classical case. We formally define $\widetilde{f}$ by taking a formal antiderivative (in the classical sense) on each component. As $1-k=1/2$, we have in fact $\widetilde{f_i}=\sqrt{\theta}f_i$. We would like to determine an alternative way to write the Eichler integral as an actual integral, so that we may use substitution and derive modularity properties of $\widetilde{f}$ from $f$. However, the integral $g(z)=c_{1/2}\int_0^{i\infty}f(z)(z-x)^{-3/2}\, dz$ no longer makes sense. To remedy this in the weight $3/2$ case, Lawrence and Zagier define another integral $f^*(x):=c_k\int_{\bar{x}}^{\infty}\frac{f(z)}{(z-x)^{\frac12}}\, dz,$ which is meaningful for $x$ in the lower half plane $\mathbb{H}^-$.
\\ 
\indent
Here we sketch their argument in the weight $1/2$ case for completeness, and as the analysis involved in our own work differs slightly. Returning to our vector-valued form $f$, recall that the definition of the Eichler integral of $f$ corresponds with $\sqrt{\theta}f$. For $x\in\mathbb{H}^-$, we define
\begin{equation}
 f^*(x)=\left(\frac{-i}{\pi(1+i)}\right)\cdot\int_{\bar{x}}^{i\infty}\frac{f(z)}{(z-x)^{\frac32}}\, dz.
 \end{equation}
To evaluate this integral, use absolute convergence to exchange the integral and the sum, and note that for $q_z=e^{2\pi i z},$
\begin{equation}
\int_{\bar{x}}^{i\infty}\frac{q_z^n}{(z-x)^{\frac32}}\, dz=\left((2+2i)\pi\sqrt{n}q_z^n\operatorname{erfi}\left((1+i)\sqrt{\pi n(z-x)}\right)-\frac{2q_x^n}{(z-x)^{\frac12}}\right)\bigg|_{z=\bar{x}}^{i\infty},
\end{equation}
where $\operatorname{erfi}(x)$ is the \emph{imaginary error function}. As in \cite{L-Zagier}, we have that $\widetilde{f}(x+i y)=f^*(x-iy)$ as full asymptotic expansions for $x\in\Q$, $0<y\in\R$. To see this, note that at the lower limit, the antiderivative vanishes as $y\rightarrow 0$ as $\operatorname{erfi}(0)=0$ and although the square root in the denominator goes to zero, for each rational at which we are evaluating our ``strange'' series, the corresponding theta functions vanish to infinite order, which makes this term converge. For the upper limit, the square root term immediately vanishes, and we use the fact that $\lim_{x\rightarrow \infty}\operatorname{erfi}(1+i)\sqrt{i x+y}=i$ for $x,y\in\R$. 
\\
\indent
Thus, as in \cite{L-Zagier}, we have that $\widetilde{f}(x)=f^*(x)$ to infinite order at rational points. In the case of $\theta_1$, we have that $\widetilde{\theta}_1(x)=\theta^*(x)$, but for $\theta_2$ and $\theta_3$ we have to divide by $4=\sqrt{16}$ due to the non-integrality of the powers of $q$ in order to agree with the definition of $\widetilde{\theta}_i$.  Using this together with Corollary 2.1, in all cases we find that $\theta_i^S(q)=\theta_i^*(q)$  at roots of unity where both sides are defined. Now, the modularity properties for the integral follow \emph{mutatis mutandis} from \cite{L-Zagier} using the modularity of $f$  and a standard $u$-substitution. More precisely, suppose $f(-1/z)(z)^{-\frac{1}{2}}=M_Sf(z)$. Then we have shown that the following modularity properties hold for $f^*(z)$ when $z\in\mathbb{H}^-$, and hence also hold for $\widetilde{f}(z)$ for each component at appropriate roots of unity where each ``strange'' function is defined. By this, we mean that the modularity conditions in the following proposition can be expressed as six equations, and each of these equations is true precisely where the corresponding ``strange'' series make sense. 
\begin{prop}
If $g(x):=\left(\frac{-i}{\pi(1+i)}\right)\cdot\int_{0}^{i\infty}\frac{f(z)}{(z-x)^{\frac32}}\, dz$, then
\[
\left(\frac{x}{-i}\right)^{-\frac{3}{2}}f(-1/x)+M_Sf(x)=M_Sg(x).
\]
\end{prop} 
It is also explained in \cite{L-Zagier} why $g_{\alpha}(z)$ is a smooth function for $\alpha\in\R$. Although $g(x)$ is \emph{a priori} only defined in $\mathbb{H}^-$, we may take any path $L$ connecting $0$ to $i\infty$. Then we can holomorphically continue $g(x)$ to all of $\C-L$. Thus, we obtain a continuation of $g$ which is smooth on $\R$ and analytic on $\R-\{0\}$.
\section{Proof of Theorem \ref{mainthm}}
Here we complete the proofs of parts (1) and (2) of Theorem \ref{mainthm}.
\subsection{Proof of Theorem \ref{mainthm} (1)}
We show that at appropriate roots of unity, our ``strange'' functions $\theta_i^S$ are reflections of $q$-series which are convergent on $\mathbb{H}$. 
Using (\ref{F9F10}), it suffices to show for $\theta_1^S$ that $\sum_{n=0}^{\infty}\frac{(q^{-1};q^{-1})_n}{(-q^{-1};q^{-1})_n}$ agrees at odd roots of unity with a $q$-series convergent when $|q|<1$. To factor out inverse powers of $q$, we observe that 
\begin{equation}\label{inversion}
(a^{-1};q^{-\alpha})_n=(-1)^na^nq^{\frac{\alpha n(n-1)}{2}}(a;q^{\alpha})_n.
\end{equation}
Applying this identity to the numerator and denominator term-by-term, we have at odd order roots of unity
\begin{equation}\theta_1^S(q^{-1})=\sum_{n=0}^{\infty}(-1)^n\frac{(q;q)_n}{(-q;q)_n}=2\displaystyle\sum_{n=0}^{\infty}\frac{q^{2n+1}(q;q)_{2n}}{(1+q^{2n+1})(-q;q)_{2n}}.
\end{equation}
The series on the right-hand side is clearly convergent for $|q|<1$, and results from pairing consecutive terms of the left-hand series as follows:
\begin{equation*}\frac{(q;q)_{2n}}{(-q;q)_{2n}}-\frac{(q;q)_{2n+1}}{(-q;q)_{2n+1}}=\frac{(q;q)_{2n}}{(-q;q)_{2n}}\left(1-\frac{1-q^{2n+1}}{1+q^{2n+1}}\right)=\frac{2q^{2n+1}(q;q)_{2n}}{(1+q^{2n+1})(-q;q)_{2n}}.
\end{equation*}
\begin{rmk}Alternatively, one can show the convergence of $\theta_1^S(q^{-1})$ by letting $a=1,b=-1,t=-1$ in Fine's identity \cite{Fine}
\begin{equation}\sum_{n=0}^{\infty}\frac{(aq;q)_n}{(bq;q)_n}(t)^n=\frac{1-b}{1-t}+\frac{b-atq}{1-t}\sum_{n=0}^{\infty}\frac{(aq;q)_n}{(bq;q)_n}(tq)^n,
\end{equation}\end{rmk}
giving \begin{equation}\label{inversefine}\theta_1^S(q^{-1})=1+\frac{q-1}{2}\sum_{n=0}^{\infty}\frac{(q;q)_n}{(-q;q)_n}(-q)^n\end{equation} which also converges for $|q|<1$. 
\\
\indent
Similarly, we use (4.1) to study $\theta_2^S,\theta_3^S$. Note that it suffices by (\ref{F9F10}) to study $\sum_{n=0}^{\infty}\frac{(q^{-2};q^{-2})_n}{(q^{-3};q^{-2})_{n}}$. Factorizing as above, we find that 
\begin{equation}
\sum_{n=0}^{\infty}\frac{(q^{-2};q^{-2})_n}{(q^{-3};q^{-2})_n}=\sum_{n=0}^{\infty}\frac{q^n(q^{2};q^{2})_n}{(q^{3};q^{2})_n}
,
\end{equation}
the right-hand side of which is clearly convergent on $\mathbb{H}$. We note that in general, similar inversion formulae result from applying (\ref{inversion}) to diverse $q$-series and other expressions involving eta functions, $q$-Pochhammer symbols and the like.
\subsection{Proof of Theorem \ref{mainthm} (2)}
\begin{proof}Here we complete the proof of Theorem \ref{mainthm} using the results of $\S2$ and $\S3$. 
Note that by the Corollary (2.1) to the sums of tails formulae of Andrews, Jim{\'e}nez-Urroz, and Ono \cite{A-J-O}, each component of $H(q)$ agrees to infinite order at rational numbers with a multiple of the corresponding Eichler integral. By the discussion of Eichler integrals in $\S3$, the value of each $\widetilde{\theta_i}$ agrees at rationals with the value of the corresponding $\theta_i^*$.  Therefore, by the discussion of the modularity properties of $\theta_i^*$, we need only to describe the modularity of $H(q)$. This is simple to check using the usual transformation laws
\begin{equation}
\eta(z+1)=\zeta_{24}\eta(z),
\end{equation}
\begin{equation}
\eta(-1/z)=\left(\frac{z}{i}\right)^{\frac12}\eta(z),
\end{equation}
and (\ref{eta_identity1}). Hence we see that 
\begin{equation}
H(z+1)=\begin{pmatrix}1&0&0\\0&0&\zeta_{12}\\0&\zeta_{24}&0 \end{pmatrix}H(z)
,
\end{equation}
\begin{equation}
H(-1/z)=\left(\frac{z}{i}\right)^{\frac12}\begin{pmatrix}0&\sqrt{2}&0\\1/\sqrt{2}&0&0\\0&0&1\end{pmatrix}H(z)
,
\end{equation}
and the corresponding transformations of $\theta_i^*$ follow.
\end{proof}
\section{Proof of Corollary 1.1}
\begin{proof}
The proof of Corollary 1.1 is a generalization of and proceeds similarly to the proofs of Theorems 4 and 5 of \cite{A-J-O}. As the sums of tails identities in Theorem 2.1 show that the ``strange'' functions $F_9$ and $F_{10}$ agree \emph{to infinite order} with the half derivatives of $F_9$ and $F_{10}$ at the roots of unity under condsideration, the coefficients in the asymptotic expansion of $H_i(t,\zeta)$ for $i=9,10$ agree up to a constant factor with the coefficients of the asymptotic expansion of $\sqrt{\theta}F_i(\zeta e^{-t})$. Recalling the classical theta series expansions for $F_i$ in (1.6), the first part of Corollary 1.1 follows immediately from the following well-known fact:
 
 \begin{lemma}[Proposition 5 of \cite{Hikami}]
 Let $\chi(n)$ be a periodic function with mean value zero and $L(s,\chi):=\sum_{n=0}^{\infty}\chi(n)n^{-s}$. As $t\searrow0$, we have
 
\[\sum_{n=0}^{\infty}n\chi(n)e^{-n^2t}\sim\sum_{n=0}^{\infty}L(-2n-1,\chi)\frac{(-t)^n}{n!}.\]
 
 \end{lemma}
The proof follows from taking a Mellin transform, making a change of variables, and picking up residues at negative integers. The assumption on the coefficients $\chi(n)$ assures that $L(s,\chi)$ can be analytically continued to $\C$. The mean value zero condition is easily checked in our case; for example for $F_9$ one needs to verify that $\{(-\zeta)^{n^2}\}_{n\geq0}$ is mean value zero for $\zeta$ a primitive order $2k+1$ root of unity, and for $F_{10}$ one must check that $\{\zeta^{\frac{(2n+1)^2}{8}}\}_{n\geq0}$ is mean value zero for an even order root of unity $\zeta$. These may both be checked using well-known results for the generalized quadratic Gauss sum
\begin{equation}
G(a,b,c):=\displaystyle\sum_{n=0}^{c-1}e\left(\frac{an^2+bn}{c}\right).
\end{equation}
In particular, for $F_9$, for an odd order root of unity $\zeta$, $-\zeta$ is primitive of order $k$ where $k\equiv2\pmod4$, so we need that $G(a,0,k)=0$ when $k\equiv2\pmod4$, which fact is well known. For $F_{10}$, we may use the standard fact that $G(a,b,c)=0$ whenever $4|c$, $(a,c)=1$, and $0<b\in2\Z+1$ to obtain our result. This Gauss sum calculation follows, for instance, by using the multiplicative property of Gauss sums together with an application of Hensel's lemma. 

In the case of $F_{10}$, note that the formula for $H_{10}(t,\alpha)$ is obtained by substituting $q=\zeta e^{-t}$ into the ``strange'' function for $F_{10}$ after letting $q\rightarrow q^{\frac18}.$ A simple change of variables in the Mellin transform in the foregoing proof of the present Lemma adjusts for the $1/8$ powers by giving an extra factor of $8^s$ before taking residues.
\end{proof}

\bibliographystyle{amsplain}
\bibliography{quantum}

\end{document}